\documentclass[11pt,reqno]{amsart}
\usepackage{amsmath}
\usepackage{amssymb}
\usepackage{amstext}
\usepackage{amsfonts}
\usepackage{amsthm}
\usepackage{ifthen}
\usepackage{xspace}
\usepackage{graphicx}
\usepackage{mathrsfs}
\usepackage{enumerate}
\usepackage{breakcites}

\numberwithin{equation}{section}  

\newtheorem{punkt}{}[section]

\theoremstyle{plain}
\newtheorem{corollary}[punkt]{Corollary}
\newtheorem{lemma}[punkt]{Lemma}
\newtheorem{proposition}[punkt]{Proposition}
\newtheorem{theorem}[punkt]{Theorem}

\theoremstyle{definition}

\theoremstyle{plain}
\newtheorem*{corollary*}{Corollary}
\newtheorem*{lemma*}{Lemma}
\newtheorem*{proposition*}{Proposition}
\newtheorem*{theorem*}{Theorem}

\theoremstyle{definition}
\newtheorem*{remark*}{Remark}
\newtheorem*{remarks*}{Remarks}
\newtheorem*{example*}{Example}
\newtheorem*{examples*}{Examples}
\newtheorem*{definition*}{Definition}
\newtheorem*{conjecture*}{Conjecture}
\newtheorem*{assumption*}{Assumption}
\newtheorem*{assumptions*}{Assumptions}
\newtheorem*{construction*}{Construction}

\def\mynat{\mathbb{N}}

\def\myreal{\mathbb{R}}

\def\mylebesgue{\lambda \mskip -8mu \lambda}

\def\myf{\mathcal{F}}

\def\myo{\mathcal{O}}

\def\ee{\mathbb{E}}
\def\pp{\mathbb{P}}

\def\re{\qopname\relax{no}{Re}\,}
\def\im{\qopname\relax{no}{Im}\,}

\def\eg{e.g.\@\xspace}
\def\ie{i.e.\@\xspace}

\def\sp#1#2{\langle{#1},{#2}\rangle}

\def\i{\qopname\relax{no}{\bf i}}

\begin{document}

\title{On the Second-Order Correlation Function \\[5pt] of the Characteristic Polynomial \\[5pt] of a Hermitian Wigner Matrix}

\author{F. G\"otze$^{1}$}
\address{Friedrich G\"otze, Fakult\"at f\"ur Mathematik, Universit\"at Bielefeld,
Postfach 100131, 33501 Bielefeld, Germany}
\email{goetze@math.uni-bielefeld.de}

\thanks{1) Supported by CRC 701 ``Spectral Structures and Topological Methods in Mathematics''}

\author{H. K\"osters}
\address{Holger K\"osters, Fakult\"at f\"ur Mathematik, Universit\"at Bielefeld,
Postfach 100131, 33501 Bielefeld, Germany}
\email{hkoesters@math.uni-bielefeld.de}

\date{December 19, 2007}

\begin{abstract}
We consider the asymptotics of the second-order correlation function
of the characteristic polynomial of a random matrix.
We show that the known result 
for a random matrix from the Gaussian Unitary Ensemble
essentially continues to hold
for a general Hermitian Wigner matrix.
Our proofs rely~on an explicit formula
for the exponential generating function
of the second-order correlation function
of the characteristic polynomial.
\end{abstract}

\maketitle

\markboth{F. G\"otze and H. K\"osters}{Characteristic Polynomials of Hermitian Wigner Matrices}

\section{Introduction}

The characteristic polynomials of random matrices 
have attracted considerable interest in the last years,
a major reason being the striking similarities between 
the (asymptotic) moments of the characteristic polynomial 
of a random matrix from the Circular Unitary Ensemble (CUE) 
and the (asymptotic) moments of the value distribution
of the Riemann zeta function along its critical line
(see \textsc{Keating} and \textsc{Snaith} \cite{KS1}).
These findings have inspired several authors to investigate 
the moments and correlation functions of the characteristic polynomial
also for other random matrix ensembles
(see \eg \textsc{Br\'ezin} and \textsc{Hikami} \cite{BH1,BH2},
\textsc{Mehta} and \textsc{Normand} \cite{MN},
\textsc{Strahov} and \textsc{Fyodorov} \cite{SF3},
\textsc{Baik}, \textsc{Deift} and \textsc{Strahov} \cite{BDS},
\textsc{Borodin} and \textsc{Strahov} \cite{BS}).

In this paper, we consider the second-order moment and correlation function
of the characteristic polynomial of a general (Hermitian) Wigner matrix:
Let $Q$ be a~probability distribution on the real line such that 
\begin{align}
\label{momentconditions}
\int x \ Q(dx) = 0 \,,
\quad
a := \int x^2 \ Q(dx) = 1/2 \,,
\quad
b := \int x^4 \ Q(dx) < \infty \,,
\end{align}
and let $(X_{ii} / \sqrt{2})_{i \in \mynat}$,
$(X^{\re}_{ij})_{i<j ,\, i,j \in \mynat}$
and
$(X^{\im}_{ij})_{i<j ,\, i,j \in \mynat}$
be independent families of independent random variables
with distribution $Q$ on some probability space $(\Omega,\myf,\pp)$.
Also,
let $X_{ij} := X^{\re}_{ij} + \i X^{\im}_{ij}$
and $X_{ji} := X^{\re}_{ij} - \i X^{\im}_{ij}$
for $i < j$, $i,j \in \mynat$.
Then, for any $N \in \mynat$, 
the (Hermitian) Wigner matrix of size $N \times N$
is~given by $X_N = (X_{ij})_{1 \leq i,j \leq N}$,
and the second-order correlation function 
of the characteristic polynomial is~given by
\begin{align}
\label{correlationfunction}
f(N;\mu,\nu) := \ee(\det (X_N - \mu I_N) \cdot \det (X_N - \nu I_N)) \,,
\end{align}
where $\mu,\nu$ are real numbers 
and $I_N$ denotes the identity matrix of size $N \times N$.
We~are interested in the asymptotics 
of the values $f(N;\mu_N,\nu_N)$ as $N \to \infty$, 
where $\mu_N,\nu_N$ depend on $N$ in some suitable fashion.

\def\fg{f_{\text{GUE}}}
\def\fu{f_{\text{UE}}}

In the special case where $Q$ is the Gaussian distribution
with mean $0$ and variance $1/2$,
the distribution of the random matrix $X_N$
is the so-called Gaussian Unitary Ensemble (GUE).
(See \eg \textsc{Forrester} \cite{Fo} or \textsc{Mehta} \cite{Me}
and the references cited therein.
However, let it be noted that these authors
use the variance $1/4$ instead of $1/2$, so that
we have to do some rescalings when using their results.)
A~remarkable feature of the GUE is that the joint distribution 
of the eigenvalues of the random matrix $X_N$ is known explicitly: 
It is given by
$$
  P_N(d\lambda_1,\hdots,d\lambda_N)
= Z_N^{-1} \cdot \prod_{1 \le i < j \le N} (\lambda_j - \lambda_i)^2 \cdot \prod_{i=1}^{N} e^{-\lambda_i^2/2} \ \mylebesgue^{N}(d\lambda_1,\hdots,d\lambda_N) \,,
$$
where $\mylebesgue^N$ denotes the $N$-dimensional Lebesgue measure
on $\myreal^N$ and $Z_N$ denotes the normalizing factor
$
Z_N = (2\pi)^{N/2} \cdot \prod_{k=1}^{N} k!\,.
$
%
%
Thus, the correlation function 
of the characteristic polynomial can be written as
$$
\fg(N;\mu,\nu) := \int_{\myreal^N} \ \prod_{i=1}^{N} (\lambda_i - \mu) \ \prod_{i=1}^{N} (\lambda_i - \nu) \ P_N(d\lambda_1,\hdots,d\lambda_N) \,,
$$
from which it follows (see \eg the proof of Proposition~4.3 in \textsc{Forrester} \cite{Fo})
that
$$
\fg(N;\mu,\nu) = \frac{\sp{p_N}{p_N}}{e^{-(\mu^2+\nu^2)/4}} \cdot K_{N+1}(\mu,\nu) \,.
$$
Here, the scalar product $\sp{\,\cdot\,}{\,\cdot\,}$ is given by 
$\sp{\varphi}{\psi} := \int_{-\infty}^{+\infty} \varphi(x) \, \psi(x) \, e^{-x^2/2} \ dx$,
the~$p_{k}$ are the monic orthogonal polynomials associated with 
this scalar product (\ie, up to scaling, the Hermite polynomials), 
and the kernel $K_N$ is given by
$$
K_N(x,y) := e^{-(x^2+y^2)/4} \sum_{k=1}^{N} \frac{p_{k-1}(x) p_{k-1}(y)}{\sp{p_{k-1}}{p_{k-1}}} \,.
$$
Using this representation, it is possible to obtain 
asymptotic approximations of the values $\fg(N;\mu_N,\nu_N)$
from the corresponding asymptotics of the Hermite poly\-nomials
(see \eg Section~8.22 in \textsc{Szeg\"o} \cite{Sz}).
It turns out that
$$
\lim_{N \to \infty} \sqrt{\frac{\pi}{2N}} \cdot \frac{1}{N!} \cdot \fg \left( N; \frac{\pi\mu}{\sqrt{N}},\frac{\pi\nu}{\sqrt{N}} \right) = \frac{\sin \pi(\mu-\nu)}{\pi(\mu-\nu)}
$$
(see \eg the proof of Proposition 4.14 in \textsc{Forrester} \cite{Fo})
and, for $\xi \in (-2,+2)$,
$$
\lim_{N \to \infty} \sqrt{\frac{1}{2 \pi N}} \cdot \frac{1}{N!} \cdot e^{-N\xi^2/2} \cdot \fg(N;\sqrt{N}\xi,\sqrt{N}\xi) = \frac{1}{2\pi} \sqrt{4-\xi^2}
$$
(see \eg the derivation of the semi-circle law 
\pagebreak[1]
in Chapter~4.3 in \textsc{Forrester} \cite{Fo}).
More generally, it is known that, for $\xi \in (-2,+2)$,
\begin{multline*}
\lim_{N \to \infty} \sqrt{\frac{1}{2 \pi N}} \cdot \frac{1}{N!} \cdot e^{-N\xi^2/2} \cdot \fg \left( N; \sqrt{N} \xi + \frac{\mu}{\sqrt{N} \varrho(\xi)} , \sqrt{N} \xi + \frac{\nu}{\sqrt{N} \varrho(\xi)} \right) \quad\\[+3pt] = e^{\xi(\mu+\nu)/2\varrho(\xi)} \cdot \varrho(\xi) \cdot \frac{\sin{\pi(\mu-\nu)}}{\pi(\mu-\nu)}
\end{multline*}
(see \eg Section~2.1 in \textsc{Strahov} and \textsc{Fyodorov} \cite{SF3}),
\pagebreak[2]
where $\varrho(\xi) := \tfrac{1}{2\pi} \sqrt{4-\xi^2}$
denotes the density of the semi-circle law.
Note that this formula includes the preceding two formulas
as special cases.

Even more, it turns out that a similar result holds for the correlation function
(of any even order $2,4,6,\hdots$)
of the characteristic polynomial of a random matrix
from the larger class of unitary-invariant ensembles
(see \eg Section~2.1 in \textsc{Strahov} and \textsc{Fyodorov} \cite{SF3}).
In this respect, it is interesting to note 
that the emergence of the sine kernel is ``universal'' in that 
it is independent of the~particular choice of the potential function
of the unitary-invariant ensemble.
In contrast to that, most of the other factors in the above result 
for the~GUE have to be replaced by potential-specific factors.

It is well-known that the GUE is a special case
not only of a unitary-invariant ensemble
but also of a (Hermitian) Wigner ensemble
as described at the beginning of this section.
The purpose of this paper is to show
that the above result for the~GUE
can also be generalized in this direction.
More precisely, our main result is as follows:

\begin{theorem}
\label{maintheorem}
Let $Q$ be a probability distribution on the real line
satisfying (\ref{momentconditions}),
let $f$ be defined as in (\ref{correlationfunction}),
let $\xi \in (-2,+2)$, and let $\mu,\nu \in \mathbb{R}$.
Then we~have
\begin{multline*}
\lim_{N \to \infty} \sqrt{\frac{1}{2 \pi N}} \cdot \frac{1}{N!} \cdot e^{-N\xi^2/2} \cdot f \left( N; \sqrt{N} \xi + \frac{\mu}{\sqrt{N}\varrho(\xi)} , \sqrt{N} \xi + \frac{\nu}{\sqrt{N}\varrho(\xi)} \right) \qquad \\ = \exp \Big( b - \tfrac{3}{4} \Big) \cdot e^{\xi(\mu+\nu) / 2\varrho(\xi)} \cdot \varrho(\xi) \cdot \frac{\sin{\pi(\mu-\nu)}}{\pi(\mu-\nu)} \,,
\end{multline*}
where $\varrho(\xi) := \tfrac{1}{2\pi} \sqrt{4-\xi^2}$ and $\sin 0 / 0 := 1$.
\end{theorem}

Specifically for the Gaussian distribution with mean $0$ and variance $1/2$,
we~have $b = \tfrac{3}{4}$, so that we re-obtain the above result
for the GUE.

Furthermore, we see that for general Wigner matrices,
the appropriately rescaled correlation function
of the characteristic polynomial asymptotically factorizes
into the universal sine kernel,
a universal factor 
involving the density of the semi-circle law,
and a non-universal factor
depending only on the fourth moment $b$,
or the fourth cumulant $b - \tfrac{3}{4}$,
of the underlying distribution $Q$.

In particular, it follows immediately that if we normalize 
the correlation function of the characteristic polynomial
by means of its second moment, we obtain the following
universality result:

\begin{corollary}
\label{maincorollary}
Under the assumptions of Theorem \ref{maintheorem}, we have
\begin{align*}
\lim_{N \to \infty} \frac{\ee \big ( D_N(\xi,\mu) \, D_N(\xi,\nu) \big)}{\sqrt{\, \ee D_N(\xi,\mu)^2 \,} \, \sqrt{\, \ee D_N(\xi,\nu)^2 \,}} 
= 
\frac{\sin{\pi(\mu-\nu)}}{\pi(\mu-\nu)} \,,
\end{align*}
where $D_N(\xi,\eta) := \det \Big( X_N - \left( \sqrt{N} \xi + \frac{\eta}{\sqrt{N} \varrho(\xi)} \right) I_N \Big)$.
\end{corollary}

Moreover, it can be shown that the asymptotics remain unchanged
if we replace the correlation function $f(N;\mu,\nu)$
by the ``true'' correlation (in the sense of prob\-ability)
of the characteristic polynomial,
\begin{align}
\label{correlationfunction2}
\widetilde{f}(N;\mu,\nu) := \ee((\det (X_N - \mu I_N) - \ee \det (X_N - \mu I_N)) \qquad\qquad\nonumber\\ \,\cdot\, (\det (X_N - \nu I_N) - \ee \det (X_N - \nu I_N))) \,.
\end{align}
We then have the following result:

\begin{proposition}
\label{maintheorem2}
Let $Q$ be a probability distribution on the real line
satisfying~(\ref{momentconditions}),
let $\widetilde{f}$ be defined as in (\ref{correlationfunction2}),
let $\xi \in (-2,+2)$, and let $\mu,\nu \in \mathbb{R}$.
Then we~have
\begin{multline*}
\lim_{N \to \infty} \sqrt{\frac{1}{2 \pi N}} \cdot \frac{1}{N!} \cdot e^{-N\xi^2/2} \cdot \widetilde{f} \left( N; \sqrt{N} \xi + \frac{\mu}{\sqrt{N}\varrho(\xi)} , \sqrt{N} \xi + \frac{\nu}{\sqrt{N}\varrho(\xi)} \right) \qquad \\ = \exp \Big( b - \tfrac{3}{4} \Big) \cdot e^{\xi(\mu+\nu) / 2\varrho(\xi)} \cdot \varrho(\xi) \cdot \frac{\sin{\pi(\mu-\nu)}}{\pi(\mu-\nu)} \,,
\end{multline*}
where $\varrho(\xi) := \tfrac{1}{2\pi} \sqrt{4-\xi^2}$ and $\sin 0 / 0 := 1$.
\end{proposition}

Similarly as before, normalizing the correlation
of the characteristic polynomial by means of its variance
leads to the following universality result
for the correlation coefficient:

\begin{corollary}
\label{corollary2}
Under the assumptions of Proposition \ref{maintheorem2}, we have
\begin{align*}
\lim_{N \to \infty} \frac{\ee \big ( \widetilde{D}_N(\xi,\mu) \, \widetilde{D}_N(\xi,\nu) \big)}{\sqrt{\, \ee \widetilde{D}_N(\xi,\mu)^2 \,} \, \sqrt{\, \ee \widetilde{D}_N(\xi,\nu)^2 \,}} 
= 
\frac{\sin{\pi(\mu-\nu)}}{\pi(\mu-\nu)} \,,
\end{align*}
where $\widetilde{D}_N(\xi,\eta) := \det \Big( X_N - \left( \sqrt{N} \xi + \frac{\eta}{\sqrt{N} \varrho(\xi)} \right) I_N \Big) - \ee \det \Big( X_N - \left( \sqrt{N} \xi + \frac{\eta}{\sqrt{N} \varrho(\xi)} \right) I_N \Big)$.
\end{corollary}

The proofs of the above-mentioned results on the correlation function
of the characteristic polynomial of a random matrix
from the GUE (or another unitary-invariant ensemble)
heavily depend on the special structure
of the joint distribution of the eigenvalues.
However, such a structure seems not to be available 
for general Wigner matrices.

Instead, we start from recursive equations
for the correlation function of the characteristic polynomial
(as well as some closely related correlation functions),
derive an explicit expression
for the associated exponential generating function
and~deduce all our asymptotic results from this expression.
For the determinant of a~real symmetric Wigner matrix,
a similar analysis was carried out by
\textsc{Zhurbenko} \cite{Zh}.

The crucial step in our analysis is to obtain an expression in closed form
for the exponential generating function of the correlation function.
Unfortunately, this approach has proven successful so far 
only for the second-order correlation function of the characteristic polynomial,
which explains why we do not have any results
for the higher-order correlation functions of the characteristic poly\-nomial.
Note however that for any distribution $Q$
with the first $2k$ moments identical to the Gaussian moments,
the correlation function of order $k$ of the characteristic polynomial
must be the same as that for the GUE.

This paper is organized as follows.
In Section~2, we start with the analysis of the recursive equations
for the correlation function of the characteristic poly\-nomial
and derive the explicit expression for its exponential generating function. 
Sections 3 and 4 are devoted to the proofs of Theorem \ref{maintheorem}
and Proposition \ref{maintheorem2}, respectively.

Throughout this paper, $K$ denotes an absolute constant 
which may change from one occurrence to the next.

\medskip

\textbf{Acknowledgement.}
We thank Mikhail Gordin for bringing the connection between 
the correlation function of the characteristic polynomial
of the GUE and the sine kernel to our attention.

\bigskip

\section{Generating Functions}

To simplify the notation, we adopt the following conventions:
\mbox{The determinant} of the ``empty'' (\ie, $0 \times 0$) matrix 
is taken to be $1$. If $A$ is an $n \times n$ matrix 
and $z$ is a real or complex number, we set $A - z := A - z I_n$, 
where $I_n$ denotes the $n \times n$ identity matrix.
Furthermore, if $A$ is an $n \times n$ matrix
and $i_1,\hdots,i_m$ and $j_1,\hdots,j_m$
are families of pairwise different indices
from the set $\{ 1,\hdots,n \}$,
we write $A^{[i_1,\hdots,i_m:j_1,\hdots,j_m]}$
for the $(n-m) \times (n-m)$-matrix obtained from $A$ 
by removing the rows indexed by $i_1,\hdots,i_m$
and the columns indexed by $j_1,\hdots,j_m$.
Thus, for any $n \times n$ matrix
$A = (a_{ij})_{1 \leq i,j \leq n}$ $(n \geq 1)$,
we have the identity
\begin{align}
\label{laplace-expansion}
\det(A) = \sum_{i,j=1}^{n-1} (-1)^{i+j-1} \, a_{i,n} \, a_{n,j} \, \det(A^{[n,i:n,j]}) + a_{n,n} \det(A^{[n:n]}) \,,
\end{align}
as follows by expanding the determinant
about the last row and the last column.
(For $n=1$, note that the big sum vanishes.)

\pagebreak[2]

Recall that we write $X_N$ for the random matrix $(X_{ij})_{1 \leq i,j \leq N}$,
where the $X_{ij}$ are the random variables introduced below (\ref{momentconditions}).
We will analyze the function
\begin{align*}
f(N;\mu,\nu) &:= \mathbb{E} \left( \det(X_N - \mu) \cdot \det(X_N - \nu) \right) \qquad\qquad & (N \geq 0) \,.
\end{align*}
To this purpose, we will also need the auxiliary functions
\begin{align*}
f_{11}(N;\mu,\nu) &:= \mathbb{E}(\det((X_N - \mu)^{[1:1]}) \cdot \det((X_N - \nu)^{[2:2]})) & (N \geq 2) \,,
\\
f_{11}^\chi(N;\mu,\nu) &:= \mathbb{E}(\det((X_N - \mu)^{[1:2]}) \cdot \det((X_N - \nu)^{[2:1]})) & (N \geq 2) \,,
\allowdisplaybreaks\\
f_{10}(N;\mu,\nu) &:= \mathbb{E}(\det(X_{N-1} - \mu) \cdot \det(X_N - \nu)) & (N \geq 1) \,,
\\
f_{01}(N;\mu,\nu) &:= \mathbb{E}(\det(X_N - \mu) \cdot \det(X_{N-1} - \nu)) & (N \geq 1) \,.
\end{align*}
Since $\mu$ and $\nu$ can be regarded as constants
for the purposes of this section, we~will only write
$f(N)$ instead of $f(N;\mu,\nu)$ in the sequel, etc.

\pagebreak[4]  

We have the following recursive equations:

\begin{lemma}
\label{fn-recursion}
\begin{align}
f(0) &= 1 \,, \nonumber\\
f(N) &= (1 + \mu \nu) \, f(N-1) + (2b + \tfrac{1}{2})(N-1) \, f(N-2)
\nonumber\\&\qquad \,+\, (N-1)(N-2) \, f_{11}  (N-1)
\nonumber\\&\qquad \,+\, (N-1)(N-2) \, f_{11}^\chi(N-1)
\nonumber\\&\qquad \,+\, \nu(N-1) \, f_{10}(N-1)
\nonumber\\&\qquad \,+\, \mu(N-1) \, f_{01}(N-1) & (N \geq 1) \,,
\label{fn-recursion-1}
\allowdisplaybreaks[2]\\[+5pt]
f_{11}(N) &= \mu \nu \, f(N-2) + (N-2) \, f(N-3)
\nonumber\\&\qquad \,+\, (N-2)(N-3) \, f_{11}  (N-2)
\nonumber\\&\qquad \,+\, \nu(N-2) \, f_{10}(N-2)
\nonumber\\&\qquad \,+\, \mu(N-2) \, f_{01}(N-2) & (N \geq 2) \,,
\label{fn-recursion-2}
\allowdisplaybreaks[2]\\[+5pt]
f_{11}^\chi(N) &= f(N-2) + (N-2) \, f(N-3)
\nonumber\\&\qquad \,+\, (N-2)(N-3) \, f_{11}^\chi(N-2) & (N \geq 2) \,,
\label{fn-recursion-3}
\allowdisplaybreaks[2]\\[+5pt]
f_{10}(N) &= -(N-1) \, f_{01}(N-1) - \nu \, f(N-1) & (N \geq 1) \,,
\label{fn-recursion-4}
\allowdisplaybreaks[2]\\[+5pt]
f_{01}(N) &= -(N-1) \, f_{10}(N-1) - \mu \, f(N-1) & (N \geq 1) \,.
\label{fn-recursion-5}
\end{align}
\end{lemma}

For the sake of clarity, note that these recursive equations
may contain some terms which have not been defined
(such as $f_{11}(N-1)$ for $N=1$), but this is not a problem 
since these terms occur in combination with the factor zero only.

\begin{proof}
We begin with the proof of (\ref{fn-recursion-1}).
For $N=0$, the result is clear.
For $N \geq 1$, we expand the determinants 
of the matrices $(X_N-\mu)$ and $(X_N-\nu)$
as in (\ref{laplace-expansion}) and use the independence
of the random variables $X_{ij} = \overline{X}_{ji}$ ($i \leq j$),
to the effect that
\begin{align*}
&\mskip24mu f(N) \\
&= \sum_{i,j=1}^{N-1} \sum_{k,l=1}^{N-1} (-1)^{i+j+k+l} \, \ee \left( X_{i,N} X_{N,j} X_{k,N} X_{N,l} \right) \cdot \ee \left( \det(X_{N-1}-\mu)^{[i:j]} \cdot \det(X_{N-1}-\nu)^{[k:l]} \right) \\
     &\quad \,+\, \sum_{i,j=1}^{N-1} (-1)^{i+j+1} \, \ee \left( X_{i,N} X_{N,j} \right) \cdot \ee \left( X_{N,N}-\nu \right) \cdot \ee \left(\det(X_{N-1}-\mu)^{[i:j]} \cdot \det(X_{N-1}-\nu) \right) \\
     &\quad \,+\, \sum_{k,l=1}^{N-1} (-1)^{k+l+1} \, \ee \left( X_{k,N} X_{N,l} \right) \cdot \ee \left( X_{N,N}-\mu \right) \cdot \ee \left(\det(X_{N-1}-\mu) \cdot \det(X_{N-1}-\nu)^{[k:l]} \right) \\
     &\quad \,+\, \ee \left( (X_{N,N}-\mu) (X_{N,N}-\nu) \right) \cdot \ee \left( \det(X_{N-1}-\mu) \cdot \det(X_{N-1}-\nu) \right) \,.
\end{align*}
Since the (complex-valued) random variables $X_{ij} = \overline{X_{ji}}$ ($i \le j$)
are independent with $\ee(X_{ij}) = 0$ ($i \leq j$) and $\ee(X_{ij}^2) = 0$ ($i < j$),
several of the expectations vanish, and the sum reduces to
\begin{align*}
&\mskip24mu f(N) \\
&= ( \ee X_{N,N}^2 + \mu \nu ) \cdot \ee \left( \det(X_{N-1} - \mu) \cdot \det(X_{N-1} - \nu) \right) \\
     &\quad \,+\, \sum_{i = j  =  k = l} \ee |X_{i,N}|^4 \cdot \ee \left( \det(X_{N-1}-\mu)^{[i:j]} \cdot \det(X_{N-1}-\nu)^{[k:l]} \right) \\
     &\quad \,+\, \sum_{i = j \ne k = l} \ee |X_{i,N}|^2 \cdot \ee |X_{k,N}|^2 \cdot \ee \left( \det(X_{N-1}-\mu)^{[i:j]} \cdot \det(X_{N-1}-\nu)^{[k:l]} \right) \\
     &\quad \,+\, \sum_{i = l \ne j = k} \ee |X_{i,N}|^2 \cdot \ee |X_{k,N}|^2 \cdot \ee \left( \det(X_{N-1}-\mu)^{[i:j]} \cdot \det(X_{N-1}-\nu)^{[k:l]} \right) \displaybreak[2] \\
     &\quad \,+\, \nu \sum_{i=j} \ee |X_{i,N}|^2 \cdot \ee \left(\det(X_{N-1}-\mu)^{[i:j]} \cdot \det(X_{N-1}-\nu) \right) \\
     &\quad \,+\, \mu \sum_{k=l} \ee |X_{k,N}|^2 \cdot \ee \left(\det(X_{N-1}-\mu) \cdot \det(X_{N-1}-\nu)^{[k:l]} \right) \,.
\end{align*}
(\ref{fn-recursion-1}) now follows
by noting that
$\ee X_{N,N}^2 = 1$,
$\ee |X_{i,N}|^2 = 1$,
$\ee |X_{i,N}|^4 = 2b + \tfrac{1}{2}$,
and by using symmetry.

To prove (\ref{fn-recursion-2}),
we apply the analogue of (\ref{laplace-expansion})
for the first row and the first column
to the matrices $(X_N-\mu)^{[1:1]}$ and $(X_N-\nu)^{[2:2]}$.
Using similar arguments as in the proof of (\ref{fn-recursion-1})
afterwards, we obtain
\begin{align*}
&\mskip24mu f_{11}(N) \\ 
    &= \sum_{i,j=3}^{N} \sum_{k,l=3}^{N} (-1)^{i+j+k+l} \, \ee \left( X_{i,2} X_{2,j} \right) \cdot \ee \left( X_{k,1} X_{1,l} \right) \cdot \ee \left( \det(X_N-\mu)^{[1,2,i:1,2,j]} \cdot \det(X_N-\nu)^{[1,2,k:1,2,l]} \right) \\
    &\quad \,+\, \sum_{i,j=3}^{N} (-1)^{i+j+1} \, \ee \left( X_{i,2} X_{2,j} \right) \cdot \ee \left( X_{1,1}-\nu \right) \cdot \ee \left(\det(X_N-\mu)^{[1,2,i:1,2,j]} \cdot \det(X_N-\nu)^{[1,2:1,2]} \right) \\
    &\quad \,+\, \sum_{k,l=3}^{N} (-1)^{k+l+1} \, \ee \left( X_{k,1} X_{1,l} \right) \cdot \ee \left( X_{2,2}-\mu \right) \cdot \ee \left(\det(X_N-\mu)^{[1,2:1,2]} \cdot \det(X_N-\nu)^{[1,2,k:1,2,l]} \right) \\
    &\quad \,+\, \ee \left( X_{2,2}-\mu \right) \cdot \ee \left( X_{1,1}-\nu \right) \cdot \ee \left( \det(X_N-\mu)^{[1,2:1,2]} \cdot \det(X_N-\nu)^{[1,2:1,2]} \right) \allowdisplaybreaks \\
    &= \mu \nu \cdot \ee \left( \det(X_N - \mu)^{[1,2:1,2]} \cdot \det(X_N - \nu)^{[1,2:1,2]} \right) \\
    &\quad \,+\, \sum_{i = j  =  k = l} \ee |X_{i,2}|^2 \cdot \ee |X_{k,1}|^2 \cdot \ee \left( \det(X_N-\mu)^{[1,2,i:1,2,j]} \cdot \det(X_N-\nu)^{[1,2,k:1,2,l]} \right) \\
    &\quad \,+\, \sum_{i = j \ne k = l} \ee |X_{i,2}|^2 \cdot \ee |X_{k,1}|^2 \cdot \ee \left( \det(X_N-\mu)^{[1,2,i:1,2,j]} \cdot \det(X_N-\nu)^{[1,2,k:1,2,l]} \right) \\
    &\quad \,+\, \nu \sum_{i=j} \ee |X_{i,2}|^2 \cdot \ee \left(\det(X_N-\mu)^{[1,2,i:1,2,j]} \cdot \det(X_N-\nu)^{[1,2:1,2]} \right) \\
    &\quad \,+\, \mu \sum_{k=l} \ee |X_{k,1}|^2 \cdot \ee \left(\det(X_N-\mu)^{[1,2:1,2]} \cdot \det(X_N-\nu)^{[1,2,k:1,2,l]} \right)
\end{align*}
and hence (\ref{fn-recursion-2}).

The proof of (\ref{fn-recursion-3})
is similar to that of (\ref{fn-recursion-2}):
\begin{align*}
&\mskip24mu f_{11}^{\chi}(N) \\
    &= \sum_{i,j=3}^{N} \sum_{k,l=3}^{N} (-1)^{i+j+k+l} \, \ee \left( X_{i,1} X_{1,l} \right) \cdot \ee \left( X_{2,j} X_{k,2} \right) \cdot \ee \left( \det(X_N-\mu)^{[1,2,i:1,2,j]} \cdot \det(X_N-\nu)^{[1,2,k:1,2,l]} \right) \\
    &\quad \,+\, \sum_{i,j=3}^{N} (-1)^{i+j+1} \, \ee \left( X_{i,1} X_{2,j} \right) \cdot \ee \left( X_{1,2} \right) \cdot \ee \left(\det(X_N-\mu)^{[1,2,i:1,2,j]} \cdot \det(X_N-\nu)^{[1,2:1,2]} \right) \\
    &\quad \,+\, \sum_{k,l=3}^{N} (-1)^{k+l+1} \, \ee \left( X_{k,2} X_{1,l} \right) \cdot \ee \left( X_{2,1} \right) \cdot \ee \left(\det(X_N-\mu)^{[1,2:1,2]} \cdot \det(X_N-\nu)^{[1,2,k:1,2,l]} \right) \\
    &\quad \,+\, \ee \left( X_{2,1} X_{1,2} \right) \cdot \ee \left( \det(X_N-\mu)^{[1,2:1,2]} \cdot \det(X_N-\nu)^{[1,2:1,2]} \right) \allowdisplaybreaks \\
    &= \ee \left( \det(X_N - \mu)^{[1,2:1,2]} \cdot \det(X_N - \nu)^{[1,2:1,2]} \right) \\
    &\quad \,+\, \sum_{i = l  =  k = j} \ee |X_{i,1}|^2 \cdot \ee |X_{k,2}|^2 \cdot \ee \left( \det(X_{N-1}-\mu)^{[1,2,i:1,2,j]} \cdot \det(X_{N-1}-\nu)^{[1,2,k:1,2,l]} \right) \\
    &\quad \,+\, \sum_{i = l \ne k = j} \ee |X_{i,1}|^2 \cdot \ee |X_{k,2}|^2 \cdot \ee \left( \det(X_{N-1}-\mu)^{[1,2,i:1,2,j]} \cdot \det(X_{N-1}-\nu)^{[1,2,k:1,2,l]} \right) \,.
\end{align*}

For the proof of (\ref{fn-recursion-4}),
we expand the determinant of the matrix $(X_N-\nu)$
as in~(\ref{laplace-expansion})
and use similar arguments as above to obtain
\begin{align*}
f_{10}(N) 
    &= \sum_{k,l=1}^{N-1} (-1)^{k+l+1} \, \ee \left( X_{k,N} X_{N,l} \right) \cdot \ee \left(\det(X_{N-1}-\mu) \cdot \det(X_{N-1}-\nu)^{[k:l]} \right) \\
    &\quad \,+\, \ee \left( X_{N,N}-\nu \right) \cdot \ee \left( \det(X_{N-1}-\mu) \cdot \det(X_{N-1}-\nu) \right) \\
    &= - \sum_{k=l} \ee |X_{k,N}|^2 \cdot \ee \left(\det(X_{N-1}-\mu) \cdot \det(X_{N-1}-\nu)^{[k:l]} \right) \\
    &\quad \,-\, \nu \cdot \ee \left( \det(X_{N-1}-\mu) \cdot \det(X_{N-1}-\nu) \right)
\end{align*}
and hence (\ref{fn-recursion-4}).

The proof of (\ref{fn-recursion-5}) is completely analogous
to that of (\ref{fn-recursion-4}).
\end{proof}

\pagebreak[2]

The interesting (although apparently rather special) phenomenon is
that the preceding recursions can be combined into a single recursion
involving only the values $f(N)$.
To shorten the notation, we put
$$
c(N) := \frac{f(N)}{N!} \qquad (N \geq 0)
$$
and \mbox{\qquad\qquad\qquad}
$$
s(N) := \sum_{\substack{k=0,\hdots,N \\ \text{$k$ even}}} c(N-k) \qquad (N \geq 0) \,.
$$
We then have the following result:

\begin{lemma}
\label{fn-superrecursion}
The values $c(N)$ satisfy the recursive equation
\begin{align}
\label{cn-formula-3}
c(0) &= 1 \,, \\
N c(N) &= c(N-1) + N \cdot c(N-2)
    \nonumber\\&\qquad \,+\, \mu \nu \cdot \big( s(N-1) + s(N-3) \big)
    \nonumber\\&\qquad \,-\, (\mu^2 + \nu^2) \cdot s(N-2)
    \nonumber\\&\qquad \,+\, (2b-\tfrac{3}{2}) \cdot \big( c(N-2) - c(N-4) \big) & (N \geq 1) \,,
\label{cn-formula-4}
\end{align}
where all terms $c(\,\cdot\,)$ and $s(\,\cdot\,)$
with a negative argument are taken to be zero.
\end{lemma}

\begin{proof}
It is immediate from Lemma \ref{fn-recursion} that
$$
f(N) = f_{11}(N+1) + f_{11}^{\chi}(N+1) + (2b-\tfrac{3}{2}) (N-1) \cdot f(N-2)
$$
for all $N \geq 1$.
Using this relation for $N-2$ instead of $N$ (where now $N \geq 3$),
we~can substitute $f_{11}(N-1) + f_{11}^{\chi}(N-1)$
on the right-hand side of (\ref{fn-recursion-1})
to~obtain
\begin{align*}
f(N) &= (1 + \mu \nu) \, f(N-1) + (2b + \tfrac{1}{2})(N-1) \, f(N-2)
    \\&\quad \,+\, (N-1)(N-2) \, \Big( f(N-2) - (2b-\tfrac{3}{2}) (N-3) f(N-4) \Big)
    \\&\quad \,+\, \nu(N-1) \, f_{10}(N-1)
    \\&\quad \,+\, \mu(N-1) \, f_{01}(N-1)
\allowdisplaybreaks[2]\\
&= (1 + \mu \nu) \, f(N-1) + N (N-1) \, f(N-2)
     \\&\quad \,+\, \nu(N-1) \, f_{10}(N-1)
     \\&\quad \,+\, \mu(N-1) \, f_{01}(N-1) 
     \\&\quad \,+\, (2b-\tfrac{3}{2}) \cdot \Big( (N-1) \, f(N-2) - (N-1)(N-2)(N-3) \, f(N-4) \Big)
\end{align*}
for all $N \geq 3$. Dividing by $(N-1)!$, it follows that
\begin{align*}
N c(N) &= (1 + \mu \nu) \cdot c(N-1) + N c(N-2)
     \\&\quad \,+\, \nu \cdot f_{10}(N-1) \,/\, (N-2)!
     \\&\quad \,+\, \mu \cdot f_{01}(N-1) \,/\, (N-2)!
     \\&\quad \,+\, (2b-\tfrac{3}{2}) \cdot \Big( c(N-2) - c(N-4) \Big)
\end{align*}
for all $N \geq 3$. \pagebreak[1] (For $N=3$, note that 
the second term in the large bracket vanishes.)
A straightforward induction 
using (\ref{fn-recursion-4}) and (\ref{fn-recursion-5}) 
shows that
\begin{align*}
f_{10}(N-1) \,/\, (N-2)! &= - \nu s(N-2) + \mu s(N-3) \,,
\\[+5pt]
f_{01}(N-1) \,/\, (N-2)! &= - \mu s(N-2) + \nu s(N-3) \,,
\end{align*}
for all $N \geq 3$, which yields the assertion for $N \geq 3$.

\pagebreak[3]

For $N < 3$, the~assertion is~verified by direct calculation,
also making use of Lemma \ref{fn-recursion}:
\begin{align*}
  c(0) = f(0) &= 1 \,.
\allowdisplaybreaks\\[+5pt]
1 c(1) = f(1) &= (1 + \mu \nu) f(0)
\\            &= (1 + \mu \nu) c(0)
\\            &= c(0) + \mu \nu s(0) \,. 
\allowdisplaybreaks\\[+5pt]
2 c(2) = f(2) &= (1 + \mu\nu) f(1) + (2b+\tfrac{1}{2}) f(0) + \nu (-\nu f(0)) + \mu (-\mu f(0))
\\            &= (1 + \mu\nu) f(1) + (2b+\tfrac{1}{2}) f(0) - (\mu^2 + \nu^2) f(0)
\\            &= (1 + \mu\nu) c(1) + (2b+\tfrac{1}{2}) c(0) - (\mu^2 + \nu^2) c(0)
\\            &= c(1) + 2 c(0) + \mu\nu c(1) - (\mu^2 + \nu^2) c(0) + (2b-\tfrac{3}{2}) c(0)
\\            &= c(1) + 2 c(0) + \mu\nu s(1) - (\mu^2 + \nu^2) s(0) + (2b-\tfrac{3}{2}) c(0) \,.
\end{align*}
\end{proof}

We can now determine the exponential generating function
of the sequence $(f(N))_{N \geq 0}$:

\begin{lemma}
\label{fn-genfun}
The exponential generating function $F(x) := \sum_{N=0}^{\infty} f(N) \, x^N \,/\, N!$ 
of the sequence $(f(N))_{N \geq 0}$ is given by
$$
F(x) = \frac{\exp \left( \mu \nu \cdot \frac{x}{1-x^2} - \tfrac{1}{2} (\mu^2 + \nu^2) \cdot \frac{x^2}{1-x^2} + b^* x^2 \right)}{(1-x)^{3/2} \cdot (1+x)^{1/2}} \,,
$$
where $b^* := b - \tfrac{3}{4}$.
\end{lemma}

\begin{proof}
It is straightforward to obtain $F(x)$
starting from (\ref{cn-formula-3}) and (\ref{cn-formula-4})
and using the basic properties of generating functions.
For the sake of completeness, we~provide a detailed proof.

To begin with, recall that $f(N) / N! = c(N)$.
Multiplying (\ref{cn-formula-4}) by $x^{N-1}$, summing over $N$
and recalling our convention concerning negative arguments, we have
\begin{align*}
  \sum_{N=1}^{\infty} N c(N) x^{N-1}
= \sum_{N=1}^{\infty} & c(N-1) x^{N-1}
+ \sum_{N=2}^{\infty} N c(N-2) x^{N-1}
\\&+ \mu \nu \left( \sum_{N=1}^{\infty} s(N-1) x^{N-1} + \sum_{N=3}^{\infty} s(N\!-\!3) x^{N-1} \right)
\\&- (\mu^2+\nu^2) \sum_{N=2}^{\infty} s(N-2) x^{N-1}
\\&+ 2b^* \left( \sum_{N=2}^{\infty} c(N-2) x^{N-1} - \sum_{N=4}^{\infty} c(N-4) x^{N-1} \right) \,,
\end{align*}
whence
\begin{multline*}
F'(x) = F(x) + (2x F(x) + x^2 F'(x)) \\ \,+\, \mu \nu \frac{1+x^2}{1-x^2} F(x) - (\mu^2+\nu^2) \frac{x}{1-x^2} F(x) + 2b^* \left( x F(x) - x^3 F(x) \right) \,.
\end{multline*}
We therefore obtain the differential equation
$$
F'(x) = \left( \frac{1 + 2x}{1-x^2} + \mu \nu \frac{1+x^2}{(1-x^2)^2} - (\mu^2+\nu^2) \frac{x}{(1-x^2)^2} + 2b^*x \right) F(x) \,,
$$
which has the solution
$$
F(x) = \frac{F_0}{(1-x)^{3/2} \cdot (1+x)^{1/2}} \exp \left( \mu \nu \frac{x}{1-x^2} - \tfrac{1}{2} (\mu^2 + \nu^2) \frac{1}{1-x^2} + b^* x^2 \right) \,.
$$
Here, $F_0$ denotes a multiplicative constant which is determined by (\ref{cn-formula-3}):
$$
F_0 = \exp(\tfrac{1}{2} (\mu^2 + \nu^2)) \,.
$$
We therefore obtain
$$
F(x) = \frac{1}{(1-x)^{3/2} \cdot (1+x)^{1/2}} \exp \left( \mu \nu \frac{x}{1-x^2} - \tfrac{1}{2} (\mu^2 + \nu^2) \frac{x^2}{1-x^2} + b^* x^2 \right) \,,
$$
which completes the proof.
\end{proof}

\bigskip

\section{The Proof of Theorem \ref{maintheorem}}

To prove Theorem \ref{maintheorem},
we will establish the following slightly more general result:

\begin{proposition}
\label{proposition-A}
Let $Q$ be a probability distribution on the real line
satisfying~(\ref{momentconditions}),
let $f$ be defined as in (\ref{correlationfunction}),
let $(\xi_N)_{N \in \mathbb{N}}$ be a sequence of real numbers
such that $\lim_{N \to \infty} \xi_N / \sqrt{N} = \xi$
for some $\xi \in (-2,+2)$,
and let $\eta \in \mathbb{C}$.
Then we~have
\begin{multline*}
\lim_{N \to \infty} \sqrt{\frac{2\pi}{N}} \cdot \frac{1}{N!} \cdot \exp(-\xi_N^2/2) \cdot f_N \left( \xi_N+\frac{\eta}{\sqrt{N}},\xi_N-\frac{\eta}{\sqrt{N}} \right)
\\
= \exp \left( b-\tfrac{3}{4} \right) \cdot \sqrt{4-\xi^2} \cdot \frac{\sin (\sqrt{4-\xi^2} \cdot \eta)}{(\sqrt{4-\xi^2} \cdot \eta)} \,,
\end{multline*}
where $\sin 0 / 0 := 1$.
\end{proposition}

It is easy to deduce Theorem \ref{maintheorem} 
from Proposition \ref{proposition-A}:

\begin{proof}[Proof of Theorem \ref{maintheorem}]
Taking
$$
\xi_N := \sqrt{N} \xi + \frac{\pi(\mu+\nu)}{\sqrt{N} \cdot \sqrt{4-\xi^2}}
\qquad\text{and}\qquad
\eta := \frac{\pi(\mu-\nu)}{\sqrt{4-\xi^2}}
$$
in Proposition \ref{proposition-A}, we have
\begin{multline*}
\lim_{N \to \infty} \sqrt{\frac{2\pi}{N}} \cdot \frac{1}{N!} \cdot \exp \left(- N \xi^2 / 2 - \pi \xi (\mu+\nu) / \sqrt{4-\xi^2} \right) \\ \,\cdot\, f_N \bigg( \sqrt{N}\xi+\frac{2\pi\mu}{\sqrt{N} \sqrt{4-\xi^2}},\sqrt{N}\xi+\frac{2\pi\nu}{\sqrt{N} \sqrt{4-\xi^2}} \bigg)
\\
= \exp \left( b-\tfrac{3}{4} \right) \cdot \sqrt{4-\xi^2} \cdot \frac{\sin \pi(\mu-\nu)}{\pi(\mu-\nu)} \,.
\end{multline*}
Multiplying by $\tfrac{1}{2\pi} \exp \left(\pi \xi (\mu+\nu) / \sqrt{4-\xi^2} \right)$
yields Theorem \ref{maintheorem}.
\end{proof}

\pagebreak[2]

It therefore remains to prove Proposition \ref{proposition-A}:

\begin{proof}[Proof of Proposition \ref{proposition-A}]
By Lemma \ref{fn-genfun}, we have
$$
\sum_{N=0}^{\infty} \frac{f(N;\mu,\nu)}{N!} z^N
=
\frac{\exp \left( \mu \nu \cdot \frac{z}{1-z^2} - \tfrac{1}{2} (\mu^2 + \nu^2) \cdot \frac{z^2}{1-z^2} + b^* z^2 \right)}{(1-z)^{3/2} \cdot (1+z)^{1/2}} \,.
$$
Thus, by Cauchy's formula, we have the integral representation
\begin{align}
\label{intrep1}
\frac{f(N;\mu,\nu)}{N!}
=
\frac{1}{2 \pi i} \int_\gamma \frac{\exp \left( \mu \nu \cdot \frac{z}{1-z^2} - \tfrac{1}{2} (\mu^2 + \nu^2) \cdot \frac{z^2}{1-z^2} + b^* z^2 \right)}{(1-z)^{3/2} \cdot (1+z)^{1/2}} \ \frac{dz}{z^{N+1}} \,,
\end{align}
where $\gamma \equiv \gamma_N$ denotes the counterclockwise circle of radius 
$R \equiv R_N = 1 - 1/N$ around the origin.
(We may and do assume that $N \geq 2$ for the rest of the~proof.) 
Setting $\mu = \xi_N + \eta / \sqrt{N}$ and $\nu = \xi_N - \eta / \sqrt{N}$,
we have
\begin{align*}
&\mskip24mu 
   \exp \left( \mu \nu \cdot \frac{z}{1-z^2} - \tfrac{1}{2} (\mu^2 + \nu^2) \cdot \frac{z^2}{1-z^2} + b^* z^2 \right) \\
&= \exp \left( (\xi_N^2 - \eta^2 / N) \cdot \frac{z}{1-z^2} - (\xi_N^2 + \eta^2 / N) \cdot \frac{z^2}{1-z^2} + b^* z^2  \right) \\
&= \exp \left( \xi_N^2 \cdot \frac{z}{1+z} - (\eta^2 / N) \cdot \frac{z}{1-z} + b^* z^2  \right) \\
&= \exp \left( \tfrac{1}{2} \xi_N^2 + \eta^2 / N \right) \cdot \exp \left( - \tfrac{1}{2} \xi_N^2 \cdot \frac{1-z}{1+z} - (\eta^2 / N) \cdot \frac{1}{1-z} + b^* z^2  \right) \,.
\end{align*}
We therefore obtain
\begin{align}
\label{intrep2}
\frac{1}{N!} \cdot f & \left( N;\xi_N+\frac{\eta}{\sqrt{N}},\xi_N-\frac{\eta}{\sqrt{N}} \right)
=
\exp \left( \tfrac{1}{2} \xi_N^2 + \eta^2 / N \right)
\nonumber\\& \,\cdot\, 
\frac{1}{2 \pi i} \int_\gamma \frac{\exp \left( - \tfrac{1}{2} \xi_N^2 \cdot \frac{1-z}{1+z} - (\eta^2 / N) \cdot \frac{1}{1-z} + b^* z^2  \right)}{(1-z)^{3/2} \cdot (1+z)^{1/2}} \ \frac{dz}{z^{N+1}} \,.
\end{align}
The idea is that the main contribution to the integral
in (\ref{intrep2}) comes from a~small neighborhood of $z=1$,
where the function
$$
h(z) := \frac{\exp \left( - \tfrac{1}{2} \xi_N^2 \cdot \frac{1-z}{1+z} - (\eta^2 / N) \cdot \frac{1}{1-z} + b^* z^2  \right)}{(1-z)^{3/2} \cdot (1+z)^{1/2}}
$$
can be well approximated by the simpler function
$$
h_0(z) := \frac{\exp(b^*)}{\sqrt{2}} \cdot \frac{\exp \left( - \tfrac{1}{4} \xi_N^2 \cdot (1-z) - (\eta^2 / N) \cdot \frac{1}{1-z} \right)}{(1-z)^{3/2}} \,.
$$
We therefore rewrite the integral in (\ref{intrep2}) as
\begin{align}
\label{intrep3}
\frac{1}{2 \pi i} \int_\gamma h(z) \ \frac{dz}{z^{N+1}}
=
I_1 + I_2 + I_3 - I_4
\end{align}
with
\begin{align}
\label{int1}
I_1 &:= \frac{1}{2 \pi i} \int_\gamma h_0(z) \ \frac{dz}{z^{N+1}} \,,
\\
\label{int2}
I_2 &:= \frac{1}{2 \pi} \int_{1/\sqrt{N}}^{2\pi-1/\sqrt{N}} h(Re^{it}) \ \frac{dt}{(Re^{it})^{N}} \,,
\\
\label{int3}
I_3 &:= \frac{1}{2 \pi} \int_{-1/\sqrt{N}}^{+1/\sqrt{N}} \left( h(Re^{it}) - h_0(Re^{it}) \right) \ \frac{dt}{(Re^{it})^{N}} \,,
\\
\label{int4}
I_4 &:= \frac{1}{2 \pi} \int_{1/\sqrt{N}}^{2\pi-1/\sqrt{N}} h_0(Re^{it}) \ \frac{dt}{(Re^{it})^{N}} \,.
\end{align}
We will show that the integral $I_1$ is the asymptotically dominant term.

\pagebreak[2]

First~of~all, note that since $\xi_N \in \myreal$, we have
\begin{align}
\label{xi-est}
\left| \exp \big( - \tfrac{1}{4} \xi_N^2 \cdot (1-z) \big) \right| = \exp \big( - \tfrac{1}{4} \xi_N^2 \cdot \re(1-z) \big) \leq 1
\end{align}
for any $z \in \mathbb{C}$ with $\re(z) \leq 1$.
Plugging in the series expansion
$$
  \exp \left( - (\eta^2 / N) \cdot \frac{1}{1-z} \right)
= \sum_{l=0}^{\infty} \frac{(-1)^l \, \eta^{2l}}{l! \, N^l} \, \frac{1}{(1-z)^l}
$$
and using uniform convergence on the contour $\gamma$
(for fixed $N \geq 2$ and $\eta \in \mathbb{C}$), we~obtain
\begin{align}
\label{I1}
\frac{I_1}{\sqrt{N}} = \frac{\exp(b^*)}{\sqrt{2}} \cdot \sum_{l=0}^{\infty} \frac{(-1)^l \, \eta^{2l}}{l! \, N^{l+1/2}} \cdot \frac{1}{2 \pi i} \int_{\gamma} \frac{\exp \left( - \tfrac{1}{4} \xi_N^2 \cdot (1-z) \right)}{(1-z)^{l+3/2}} \ \frac{dz}{z^{N+1}} \,.
\end{align}
We will show that for each $l=0,1,2,3,\hdots,$
\begin{align}
\lim_{N \to \infty}
\bigg( \frac{(-1)^l \, \eta^{2l}}{l! \, N^{l+1/2}} \cdot \frac{1}{2 \pi i} \int_{\gamma} \frac{\exp \left( - \tfrac{1}{4} \xi_N^2 \cdot (1-z) \right)}{(1-z)^{l+3/2}} \ \frac{dz}{z^{N+1}} \bigg)
\qquad\qquad\qquad\nonumber\\
=
\frac{1}{\sqrt{\pi}} \cdot \frac{(-1)^l \, \eta^{2l}}{(2l+1)!} \cdot (4-\xi^2)^{l+1/2} \,.
\label{term-by-term}
\end{align}

To begin with,
\begin{align}
&\frac{1}{2 \pi i} \int_{\gamma} \frac{\exp \left( - \tfrac{1}{4} \xi_N^2 \cdot (1-z) \right)}{(1-z)^{l+3/2}} \ \frac{dz}{z^{N+1}}
\nonumber\\&\qquad\qquad\qquad\qquad=
\frac{1}{2 \pi i} \int_{(1-1/N)-\i\infty}^{(1-1/N)+\i\infty} \frac{\exp \left( - \tfrac{1}{4} \xi_N^2 \cdot (1-z) \right)}{(1-z)^{l+3/2}} \ \frac{dz}{z^{N+1}} \,.
\label{deform}
\end{align}
In fact, for any $R' > 1$, we can replace
the contour $\gamma$ by the contour $\delta$
which consists of the line segment between the points
$R - \i \sqrt{(R')^2 - R^2}$ and $R + \i \sqrt{(R')^2 - R^2}$,
and the arc of radius $R'$ around the origin to the left of this line segment
(see~Figure~\ref{newcontour}).
\begin{figure}
\begin{center}
\includegraphics[width=6cm]{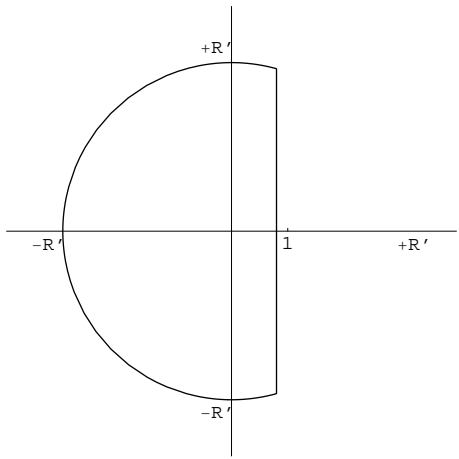}
\end{center}
\caption{The contour $\delta$.}
\label{newcontour}
\end{figure}
Now, it is easy to see that the integral along this arc is~bounded above by
$$
\frac{1}{2\pi} \cdot 2\pi R' \cdot \frac{1}{(R'-1)^{l+3/2}} \cdot \frac{1}{(R')^{N+1}}
$$
and therefore tends to zero as $R' \to \infty$, whence (\ref{deform}).

Next, performing a change of variables, we find that 
the right-hand side in (\ref{deform}) is equal to
$$
N^{l+1/2} \cdot \frac{1}{2 \pi} \int_{-\infty}^{+\infty} \frac{\exp \left( - \tfrac{1}{4} (\xi_N^2/N) \cdot (1-iu) \right)}{(1-iu)^{l+3/2}} \ \frac{du}{(1-\frac{1-iu}{N})^{N+1}} \,.
$$
Since $\lim_{N \to \infty} \xi_N / \sqrt{N} = \xi$,
it follows by the dominated convergence theorem that
\begin{multline*}
\lim_{N \to \infty}
\frac{1}{2 \pi} \int_{-\infty}^{+\infty} \frac{\exp \left( - \tfrac{1}{4} (\xi_N^2/N) \cdot (1-iu) \right)}{(1-iu)^{l+3/2}} \ \frac{du}{(1-\frac{1-iu}{N})^{N+1}}
\\ =
\frac{1}{2 \pi} \int_{-\infty}^{+\infty} \frac{\exp \left( (1 - \tfrac{1}{4} \xi^2) \cdot (1-iu) \right)}{(1-iu)^{l+3/2}} \ du \,,
\end{multline*}
which is equal to
$$
  \frac{(1-\tfrac{1}{4}\xi^2)^{l+1/2}}{\Gamma(l+3/2)}
= \frac{1}{\sqrt{\pi}} \cdot \frac{l!}{(2l+1)!} \cdot (4-\xi^2)^{l+1/2}
$$
by Laplace inversion (see \eg Chapter~24 in \textsc{Doetsch} \cite{Do}) 
and the functional equation of the Gamma function.
This proves (\ref{term-by-term}).

\pagebreak[2]

Let $\varepsilon > 0$ denote a constant such that
$\cos t \leq 1 - \varepsilon^2 t^2$ for $-\pi \leq t \leq +\pi$.
Then, for any $\alpha > 1$, we have the estimate
\begin{align}
      \int_{-\pi}^{+\pi} \frac{1}{|1-Re^{it}|^{\alpha}} \ dt
&=    \int_{-\pi}^{+\pi} \frac{1}{(1+R^2-2R \cos t)^{\alpha/2}} \ dt \nonumber\\
&\leq \int_{-\pi}^{+\pi} \frac{1}{((1-R)^2+\varepsilon^2 t^2)^{\alpha/2}} \ dt \nonumber\\
&=    N^\alpha \int_{-\pi}^{+\pi} \frac{1}{(1+N^2 \varepsilon^2 t^2)^{\alpha/2}} \ dt \nonumber\\
&=    KN^{\alpha-1} \int_{-N\varepsilon\pi}^{+N\varepsilon\pi} \frac{1}{(1+u^2)^{\alpha/2}} \ du \nonumber\\
&\leq KN^{\alpha-1} \left( 1 + \int_{1}^{\infty} \frac{1}{u^\alpha} \ du \right) \nonumber\\
&\leq KN^{\alpha-1} \left( 1 + \frac{1}{\alpha-1} \right) \,,
\label{theestimate}
\end{align}
where $K$ denotes some absolute constant which may change from~line to~line.
We~therefore obtain the bound
\begin{align*}
\sum_{l=0}^{\infty} & \left| \frac{(-1)^l \, \eta^{2l}}{l! \, N^{l+1/2}} \cdot \frac{1}{2 \pi i} \int_{\gamma} \frac{\exp \left( - \tfrac{1}{4} \xi^2_N \cdot (1-z) \right)}{(1-z)^{l+3/2}} \ \frac{dz}{z^{N+1}} \right|
\\&\leq \sum_{l=0}^{\infty} \frac{|\eta|^{2l}}{l!} \cdot \frac{1}{N^{l+1/2}} \cdot \frac{1}{2\pi} \int_{-\pi}^{+\pi} \frac{1}{|1-Re^{it}|^{l+3/2}} \ \frac{dt}{|Re^{it}|^N} 
\\&\leq K \cdot \sum_{l=0}^{\infty} \frac{|\eta|^{2l}}{l!} \cdot \left( 1 + \frac{1}{l+1/2} \right) < \infty \,,
\end{align*}
uniformly in $N \geq 2$. 
Thus, the term-by-term convergence established in (\ref{term-by-term})
entails the convergence of the complete series in (\ref{I1}),
and we~obtain
\begin{align*}
   \lim_{N\to\infty} \frac{I_1}{\sqrt{N}} 
&= \sqrt{\frac{1}{2\pi}} \cdot \exp \left( b - \tfrac{3}{4} \right) \cdot \sqrt{4-\xi^2} \cdot \sum_{l=0}^{\infty} \frac{(-1)^l \big( \sqrt{4-\xi^2} \cdot \eta \big)^{2l}}{(2l+1)!} \\
&= \sqrt{\frac{1}{2\pi}} \cdot \exp \left( b - \tfrac{3}{4} \right) \cdot \sqrt{4-\xi^2} \cdot \frac{\sin(\sqrt{4-\xi^2} \cdot \eta)}{(\sqrt{4-\xi^2} \cdot \eta)} \,.
\end{align*}
Hence, in view of (\ref{intrep2}) and (\ref{intrep3}),
the proof of Proposition \ref{proposition-A} will be complete 
once we have shown that the integrals $I_2$, $I_3$, $I_4$ 
are asymptotically negligible in the sense 
that they are of~order $o(\sqrt{N})$.

\pagebreak[2]

For the integral $I_2$, we use the estimates
($R = 1 - 1/N$, $t \in \mathbb{R}$)
\begin{align}
   \left| \exp \left( - \tfrac{1}{2} \xi_N^2 \cdot \frac{1-Re^{it}}{1+Re^{it}} \right) \right|
&= \exp \left( - \tfrac{1}{2} \xi_N^2 \cdot \re \left( \frac{1-Re^{it}}{1+Re^{it}} \right) \right) \nonumber\\
&= \exp \left( - \tfrac{1}{2} \xi_N^2 \cdot \frac{1-R^2}{1+R^2+2R\cos t} \right)
\leq 1 \,,
\label{est4}
\\[+7pt]
      \left| \exp \left( - (\eta^2 / N) \cdot \frac{1}{1-Re^{it}} \right) \right|
&\leq \exp \left( (|\eta|^2 / N) \cdot \frac{1}{|1-Re^{it}|} \right) \nonumber\\
&\leq \exp \left( (|\eta|^2 / N) \cdot \frac{1}{1-R} \right)
   =  \exp(|\eta|^2) \,,
\label{est5}
\\[+7pt]
\label{est6}
      \left| \exp \Big( b^* (Re^{it})^2 \Big) \right|
&\leq \exp \Big( |b^*| |Re^{it}|^2 \Big) 
 \leq \exp (|b^*|) \,,
\end{align}
\begin{align}
\label{est1}
|1 + Re^{it}| \geq 1 + R \cos t \geq 1 \quad\text{for}\quad \cos t \geq 0 \,,
\\[+7pt]
\label{est2}
|1 - Re^{it}| \geq 1 - R \cos t \geq 1 \quad\text{for}\quad \cos t \leq 0 \,.
\end{align}
Using these estimates, it follows that
\begin{align*}
|I_2|
&\leq 
   \frac{1}{2\pi} \int_{1/\sqrt{N}}^{\pi/2} \frac{\exp(|\eta|^2 + |b^*|)}{|1-Re^{it}|^{3/2}} \ \frac{dt}{R^N} \\
&+ \frac{1}{2\pi} \int_{\pi/2}^{3\pi/2} \frac{\exp(|\eta|^2 + |b^*|)}{|1+Re^{it}|^{1/2}} \ \frac{dt}{R^N} \\
&+ \frac{1}{2\pi} \int_{3\pi/2}^{2\pi-1/\sqrt{N}} \frac{\exp(|\eta|^2 + |b^*|)}{|1-Re^{it}|^{3/2}} \ \frac{dt}{R^N} \,.
\end{align*}
Similarly as in (\ref{theestimate}),
we have the estimates
\begin{align*}
      \int_{1/\sqrt{N}}^{\pi/2} \frac{1}{|1-Re^{it}|^{3/2}} \ dt
&\leq KN^{1/2} \int_{\sqrt{N}\varepsilon}^{\infty} \frac{1}{u^{3/2}} \ du
 \leq KN^{1/4} \,, \\[+5pt]
      \int_{\pi/2}^{3\pi/2} \frac{1}{|1+Re^{it}|^{1/2}} \ dt
&\leq KN^{-1/2} \left( 1 + \int_{1}^{N\varepsilon\pi/2} \frac{1}{u^{1/2}} \ du \right)
 \leq K \,, \\[+5pt]
      \int_{3\pi/2}^{2\pi-1/\sqrt{N}} \frac{1}{|1-Re^{it}|^{3/2}} \ dt
&\leq KN^{1/2} \int_{\sqrt{N}\varepsilon}^{\infty} \frac{1}{u^{3/2}} \ du
 \leq KN^{1/4} \,.
\end{align*}
(Recall our convention that the constant $K$ may change 
from one occurrence to the~next.) It follows that
$$
|I_2| \leq K \exp(|\eta|^2 + |b^*|) \, N^{1/4} = o(\sqrt{N}) \,.
$$

\pagebreak[2]

For the integral $I_3$, we write
$$
h(z) - h_0(z) =
\frac{\exp \left( - \tfrac{1}{4} \xi_N^2 \cdot (1-z) - (\eta^2 / N) \cdot \frac{1}{1-z} \right)}{(1-z)^{3/2}}
\cdot
\left( \tilde{h}(z) - \tilde{h}(1) \right) \,,
$$
where
$$
\tilde{h}(z) = \frac{\exp \left( - \tfrac{1}{4} \xi_N^2 \cdot \frac{(1-z)^2}{1+z} + b^* z^2 \right)}{(1+z)^{1/2}} \,,
$$
so that
\begin{multline*}
\tilde{h}'(z) = \bigg(\frac{\frac{1}{2} \xi_N^2 \cdot \frac{1-z}{1+z} + \frac{1}{4} \xi_N^2 \cdot \frac{(1-z)^2}{(1+z)^2} + 2 b^* z}{(1+z)^{1/2}} - \frac{1/2}{(1+z)^{3/2}} \bigg) \\ \,\cdot\, \exp \left( - \tfrac{1}{4} \xi_N^2 \cdot \frac{(1-z)^2}{1+z} + b^* z^2 \right) \,.
\end{multline*}
Let
$$
Z := \big\{ z \in \mathbb{C} \,|\, z=re^{i\varphi}, 1-1/N \leq r \leq 1, |\varphi| \leq 1/\sqrt{N} \big\}
$$
and note that for $z \in Z$, we have the estimate
\begin{align}
\label{est7}
\re \left( - \frac{(1-z)^2}{1+z} \right) \leq 4/N \,.
\end{align}
Indeed, with $z = re^{i\varphi}$, a simple calculation yields
$$
  \re \left( - \frac{(1-z)^2}{1+z} \right)
= \frac{- 1 + r \cos\varphi + 2r^2 - r^2 \cos 2\varphi - r^3 \cos\varphi}{1 + r^2 + 2r \cos\varphi} \,,
$$
where the denominator is clearly larger than $1$
and the numerator is bounded above~by
$$
     r \cos\varphi - r^3 \cos\varphi + r^2 - r^2 \cos 2\varphi \\
\leq r (1-r^2) \cos\varphi + r^2 - r^2 (1 - 2\varphi^2) \\
\leq 4/N \,.
$$
Since for $z=Re^{it}$ with $|t| \leq 1/\sqrt{N}$,
the line segment between the points $z$~and~$1$
is contained in the set $Z$,
it follows that
\begin{align*}
& \left| \tilde{h}(z) - \tilde{h}(1) \right|
\leq |z-1| \sup_{\alpha \in [0;1]} \left| \tilde{h}'((1-\alpha)z + \alpha) \right|
\leq |z-1| \sup_{\zeta \in Z} \left| \tilde{h}'(\zeta) \right| \\
&\quad\leq K |z-1| \, \sup_{\zeta \in Z} \, \left\{ \Big( \xi_N^2 |1-\zeta| + |b^*| + 1 \Big) \, \exp \left( \tfrac{1}{4} \xi_N^2 \cdot \re \bigg( - \frac{(1-\zeta)^2}{1+\zeta} \bigg) + |b^*| \right) \right\} \\
&\quad\leq K |z-1| \left( \xi_N^2/\sqrt{N} + |b^*| + 1 \right) \, \exp \Big( \xi_N^2 / N + |b^*| \Big) \\
&\quad\leq K(b^*,\xi^*) \, \sqrt{N} \, |z-1| \,,
\end{align*}
where the last step uses that 
$\lim_{N \to \infty} \xi_N / \sqrt{N} = \xi$,
and $K(b^*,\xi^*)$ denotes some constant depending only 
on $b^*$ and $\xi^* := (\xi_N)_{N \in \mathbb{N}}$.
Using (\ref{xi-est}), (\ref{est5}) as well as
a~similar estimate as in (\ref{theestimate}), 
we therefore obtain
\begin{align*}
      |I_3|
&\leq \frac{1}{2\pi} \int_{-1/\sqrt{N}}^{+1/\sqrt{N}} \frac{\exp(|\eta|^2)}{|1-Re^{it}|^{3/2}} \cdot \left| \tilde{h}(Re^{it}) - \tilde{h}(1) \right| \frac{dt}{R^N} \\
&\leq K(b^*,\xi^*,\eta) \, \sqrt{N} \int_{-1/\sqrt{N}}^{+1/\sqrt{N}} \frac{1}{|1-Re^{it}|^{1/2}} \ dt \\
&\leq K(b^*,\xi^*,\eta) \, \left( 1 + \int_{1}^{\sqrt{N}\varepsilon} \frac{1}{u^{1/2}} \ du \right) \\
&\leq K(b^*,\xi^*,\eta) \, N^{1/4}
   =  o(\sqrt{N}) \,,
\end{align*}
where $K(b^*,\xi^*,\eta)$ denotes some constant 
which depends only on $b^*$, $\xi^*$, and $\eta$ 
\linebreak[2] (and which may change from line to line as usual).

\pagebreak[2]

For the integral $I_4$, we can use (\ref{xi-est}), (\ref{est5}) 
as well as a similar estimate as in~(\ref{theestimate})
to obtain
\begin{align*}
      |I_4| 
&\leq \frac{1}{2\pi} \int_{1/\sqrt{N}}^{2\pi-1/\sqrt{N}} \frac{\exp(|\eta|^2 + |b^*|)}{|1-Re^{it}|^{3/2}} \ \frac{dt}{R^N} \\
&\leq K \exp(|\eta|^2 + |b^*|) \, \Biggl( \int_{1/\sqrt{N}}^{\pi} \frac{1}{|1-Re^{it}|^{3/2}} \ dt + \int_{\pi}^{2\pi-1/\sqrt{N}} \frac{1}{|1-Re^{it}|^{3/2}} \ dt \Biggr) \\
&\leq K \exp(|\eta|^2 + |b^*|) \, \Biggl( N^{1/2} \int_{\sqrt{N}\varepsilon}^{\infty} \frac{1}{u^{3/2}} \ du + N^{1/2} \int_{\sqrt{N}\varepsilon}^{\infty} \frac{1}{u^{3/2}} \ du \Biggr) \\
&\leq K \exp(|\eta|^2 + |b^*|) \, N^{1/4} = o(\sqrt{N}) \,.
\end{align*}

This completes the proof of Proposition \ref{proposition-A}.
\end{proof}

\bigskip

\section{The Proof of Proposition \ref{maintheorem2}}

Let $f(N;\mu,\nu)$ and $\widetilde{f}(N;\mu,\nu)$ be defined
as in (\ref{correlationfunction}) and (\ref{correlationfunction2}),
respectively, and note that
\begin{align}
&\mskip24mu  \widetilde{f}(N;\mu,\nu) \nonumber\\
&= \ee \big( (\det (X_N - \mu) - \ee \det (X_N - \mu)) \cdot (\det (X_N - \nu) - \ee \det (X_N - \nu)) \big) \qquad \nonumber\\
&= \ee \big( \det (X_N - \mu) \cdot \det (X_N - \nu) \big) - \ee \det (X_N - \mu) \cdot \ee \det (X_N - \nu) \nonumber\\
&= f(N;\mu,\nu) - g(N;\mu) \, g(N;\nu) \,,
\label{covariance}
\end{align}
where
$$
g(N;\mu) := \ee \det(X_N-\mu)
$$
for any $\mu \in \mathbb{R}$.
We will deduce Proposition \ref{maintheorem2} from Theorem \ref{maintheorem}
by showing that 
$f(N;\mu,\nu)$ is~asymptotically much larger than $g(N;\mu) \, g(N;\nu)$.

To this end, we need some more information about the values $g(N;\mu)$.
Similarly as in Section~2, we have the following recursive equation:

\begin{lemma}
\label{gn-recursion}
$$
g(0;\mu) = 1, \quad g(N;\mu) = - \mu g(N-1;\mu) - (N-1) g(N-2;\mu) \quad (N \geq 1) \,.
$$
\end{lemma}

\begin{proof}
For $N=0$, the claim follows from our convention
that the determinant of the empty matrix is $1$.
For $N \geq 1$, we expand the determinant of the matrix $(X_N-\mu)$ 
as in (\ref{laplace-expansion}) and use independence and symmetry 
to get
\begin{align*}
   g(N;\mu) 
&= \sum_{i,j=1}^{N-1} (-1)^{i+j+1} \, \ee \left( X_{i,N} X_{N,j} \right) \cdot \ee \left( \det(X_{N-1}-\mu)^{[i:j]} \right) \\
&\qquad \,+\, \ee \left( X_{NN}-\mu \right) \cdot \ee \left( \det(X_{N-1}-\mu) \right) \\
&= \sum_{i=j} (-1) \, \ee | X_{i,N} |^2 \cdot \ee \left( \det(X_{N-1}-\mu)^{[i:i]} \right) \\
&\qquad \,+\, \ee \left( X_{NN}-\mu \right) \cdot \ee \left( \det(X_{N-1}-\mu) \right) \\
&= - (N-1) g(N-2;\mu) - \mu g(N-1;\mu) \,. 
\end{align*}
This completes the proof of Lemma \ref{gn-recursion}.
\end{proof}

It follows from Lemma \ref{gn-recursion} that
the polynomials $g(N;\mu)$ coincide, up to scaling,
with the Hermite polynomials $H_N(x)$
(see \eg Section~5.5 in \textsc{Szeg\"o} \cite{Sz}),
which satisfy the recursive equation
$$
H_0(x) = 1, \quad H_N(x) = 2x H_{N-1}(x) - 2(N-1)H_{N-2}(x) \quad (N \geq 1) \,.
$$
(Specifically for the GUE, this is well-known,
see \eg Chapter~4 in \textsc{Forrester}~\cite{Fo}.)
The precise relationship is as follows:

\begin{lemma}
\label{gn-formula}
For any $N = 0,1,2,3,\hdots,$
$$
g(N;\mu) = (-1)^N \, 2^{-N/2} \, H_N(\mu/\sqrt{2}) \,.
$$
\end{lemma}

\begin{proof}
This follows from the recursive equations for $g(N;\mu)$ and $H_N(x)$
by a~straightforward induction on $N$.
\end{proof}

Due to Lemma \ref{gn-formula}, it is easy
to obtain the asymptotics of the values \linebreak
$g(N;\sqrt{N} \xi + \mu / \sqrt{N} \varrho(\xi))$
from the corresponding asymptotics of the Hermite poly\-nomials
(see \eg Section~8.22 in \textsc{Szeg\"o} \cite{Sz}).
For our purposes, the following estimate will be sufficient:

\begin{lemma}
\label{gn-estimate}
For $\xi \in (-2,+2)$, $\mu \in \mathbb{R}$ fixed,
$$
\left| e^{-N\xi^2/4} \, g \left( N;\sqrt{N}\xi+\frac{\mu}{\sqrt{N} \varrho(\xi)} \right) \right|
\leq
K(\xi,\mu) \, N^{-1/4} \, N!^{1/2} \,,
$$
where $K(\xi,\mu)$ is some constant depending only on $\xi$ and $\mu$.
\end{lemma}

\begin{proof}
By Theorem~8.22.9\,(a) in \textsc{Szeg\"o} \cite{Sz},
we have, for $x = \sqrt{2N+1} \, \cos \varphi$,
\begin{align*}
   e^{-x^2/2} \, H_N(x)
&= 2^{(N/2)+(1/4)} (N!)^{1/2} (N\pi)^{-1/4} (\sin \varphi)^{-1/2} \\
&\quad \,\cdot\, \left( \sin \left( ((2N+1)/4) \cdot (\sin 2\varphi - 2\varphi) + 3\pi/4 \right) + \myo(N^{-1}) \right) \,,
\end{align*}
where the $\myo$-bound holds uniformly
in $\varphi \in [\varepsilon,\pi-\varepsilon]$,
for any $\varepsilon > 0$.

\pagebreak[2]

Combining this result with Lemma \ref{gn-formula},
we obtain, for $N$ sufficiently large,
\begin{align*}
   e^{-(\sqrt{N}\xi+\mu/\sqrt{N}\varrho(\xi))^2/4} \mskip-60mu&\,\mskip+60mu g(N;\sqrt{N}\xi+\mu/\sqrt{N}\varrho(\xi)) \\
&= (-1)^N \, 2^{1/4} \, (N!)^{1/2} \, (N\pi)^{-1/4} \, (\sin \varphi)^{-1/2} \\
&\quad \,\cdot\, \left( \sin \left( ((2N+1)/4) \cdot (\sin 2\varphi - 2\varphi) + 3\pi/4 \right) + \myo(N^{-1}) \right) \,,
\end{align*}
where
$$
\varphi \equiv \varphi_N := \arccos \left( (\sqrt{N}\xi+\mu/\sqrt{N}\varrho(\xi)) \,/\, \sqrt{4N+2}\, \right)
$$
is contained in an interval of the form $[-\varepsilon,\pi-\varepsilon]$
with $\varepsilon > 0$. From this, Lemma~\ref{gn-estimate} easily follows.
\end{proof}

After these preparations we can turn to the proof of Proposition \ref{maintheorem2}:

\begin{proof}[Proof of Proposition \ref{maintheorem2}]
By Equation (\ref{covariance}) and Lemma \ref{gn-estimate},
the difference between the left-hand sides
in Theorem \ref{maintheorem} and Proposition \ref{maintheorem2}
is bounded by
\begin{multline*}
\left| \sqrt{\frac{1}{2 \pi N}} \cdot \frac{1}{N!} \cdot e^{-N\xi^2/2} \cdot g \! \left( N;\sqrt{N} \xi + \frac{\mu}{\sqrt{N}\varrho(\xi)} \right) g \! \left(N;\sqrt{N} \xi + \frac{\nu}{\sqrt{N}\varrho(\xi)} \right) \right|
\\
\leq K(\xi,\mu,\nu) \, N^{-1/2} \, N!^{-1} \, \Big( N^{-1/4} \, N!^{1/2} \Big)^2 = K(\xi,\mu,\nu) \, N^{-1} \,,
\end{multline*}
where $K(\xi,\mu,\nu)$ is some constant depending only on $\xi$, $\mu$ and $\nu$.
Thus, Proposition~\ref{maintheorem2} follows from Theorem~\ref{maintheorem}.
\end{proof}

%
%

\bigskip

\end{document}